\documentclass{article}

\usepackage{scml2024}

\usepackage{amsmath}

\usepackage[utf8]{inputenc} 
\usepackage[T1]{fontenc}    
\usepackage{hyperref}       
\usepackage{url}            
\usepackage{booktabs}       
\usepackage{amsfonts}       
\usepackage{nicefrac}       
\usepackage{microtype}      
\usepackage{cleveref}       
\usepackage{lipsum}         
\usepackage{graphicx}
\usepackage{natbib}
\usepackage{doi}
\usepackage{wrapfig}

\usepackage{amsfonts,amssymb}
\usepackage{bbm}
\usepackage{natbib}
\usepackage{amsthm}
\usepackage{algpseudocode}
\usepackage{algorithm}
\usepackage{enumerate}

\newtheorem{theorem}{Theorem}[section]
\newtheorem{lemma}[theorem]{Lemma}

\usepackage{caption}
\usepackage{subcaption}

\title{Neural Networks are Integrable}
\renewcommand{\shorttitle}{\textit{SCML2024} Template}

\date{}

\newif\ifuniqueAffiliation
\uniqueAffiliationtrue

\ifuniqueAffiliation 
\author{%
    Yucong Liu \\
    School of Mathematics \\
    Georgia Institute of Technology \\
    Atlanta, GA 30332 \\
    \texttt{yucongliu@gatech.edu} \\
}
\else
\usepackage{authblk}

\author[1]{%
	David S.~Hippocampus\thanks{\texttt{hippo@cs.cranberry-lemon.edu}}}%
\author[1,2]{%
	Elias D.~Striatum\thanks{\texttt{stariate@ee.mount-sheikh.edu}}}%
\affil[1]{Department of Computer Science, Cranberry-Lemon University, Pittsburgh, PA 15213}
\affil[2]{Department of Electrical Engineering, Mount-Sheikh University, Santa Narimana, Levand}
\fi

\hypersetup{
pdftitle={A template for the arxiv style},
pdfsubject={q-bio.NC, q-bio.QM},
pdfauthor={David S.~Hippocampus, Elias D.~Striatum},
pdfkeywords={First keyword, Second keyword, More},
}

\begin{document}
\maketitle
\begin{abstract}
In this study, we explore the integration of Neural Networks, a powerful class of functions known for their exceptional approximation capabilities. Our primary emphasis is on the integration of multi-layer Neural Networks, a challenging task within this domain. To tackle this challenge, we introduce a novel numerical method that consist of a forward algorithm and a corrective procedure. Our experimental results demonstrate the accuracy achieved through our integration approach.
\end{abstract}


\section{Introduction}
Deep learning models have demonstrated incredible power in various fields such as image and speech recognition, natural language processing, and autonomous driving in recent years. This work primarily centers on the Deep Neural Network, also known as multi-layer Perceptron or Feed-Forward Neural Network.

We begin by providing the definition of Neural Networks.
A k-layer Neural Network $\psi$ from $\mathbb{R}^{n_{0}}$ to $\mathbb{R}^{n_{k}}$ is a layer-wise structure. The $i$-the layer, for $i \in \{1,\cdots, k\}$, is defined as :
\begin{equation}
    \begin{aligned}
        y^{(i)} = W^{(i)}x^{(i-1)} + b^{(i)}, \quad
        x^{(i)} = \sigma(y^{(i)}),
    \end{aligned}
\end{equation}
and $\psi(x) = W^{(k+1)}x^{(k)} + b^{(k+1)}$.
In each layer, $W^{(i)} \in \mathbb{R}^{n_{i}\times n_{i-1}}$ is the weight matrix and $b^{(i)} \in \mathbb{R}^{n_{i}}$ is the bias vector. The input of $(i+1)$-th layer $x^{(i)} \in \mathbb{R}^{n_{i}}$ is the output of $i$-th layer, and $x^{(0)} = x$. Activation function $\sigma$ is a nonlinear point-wise function, i.e $(\sigma(y^{(i)}))_{j} = \sigma(y^{(i)}_{j})$.
One of the most well-known activation function is rectified linear unit (ReLU), which is defined as $\text{ReLU}(x) = \max(x,0)$. 

This work focuses on the topic of integration for both shallow and deep Neural Networks. Since gradient descent technique has been well used when optimizing a Neural Network, it's well studied in the literature. However the integrability of Neural Networks has garnered comparatively less attention. We aim to bridge this gap by providing explicit forms of integration for one-layer Neural Networks with any integrable activation function and deriving a piece-wise structure of the integration for multi-layer Neural Networks with ReLU activation function, along with a proposed algorithm with a corrector. 

We will introduce our basic motivation in Section \ref{2}. For the integration of ReLU Neural Networks, we separate it into two cases, one-layer case in Section \ref{3.1} and multi-layer case in Section \ref{3.2}. Our algorithm also works for Convolutional Neural Networks \citep{krizhevsky2017imagenet} and Residual Neural Networks \citep{he2016deep}, which will be introduced in Section \ref{3.3}. We will discuss about future work in Section \ref{5} and present our experiments detail in Appendix \ref{Appd:C}.

\begin{wrapfigure}[30]{R}{0.5\textwidth}
    \centering
    \begin{subfigure}[t]{0.4\textwidth}
         \centering
         \includegraphics[width=\textwidth]{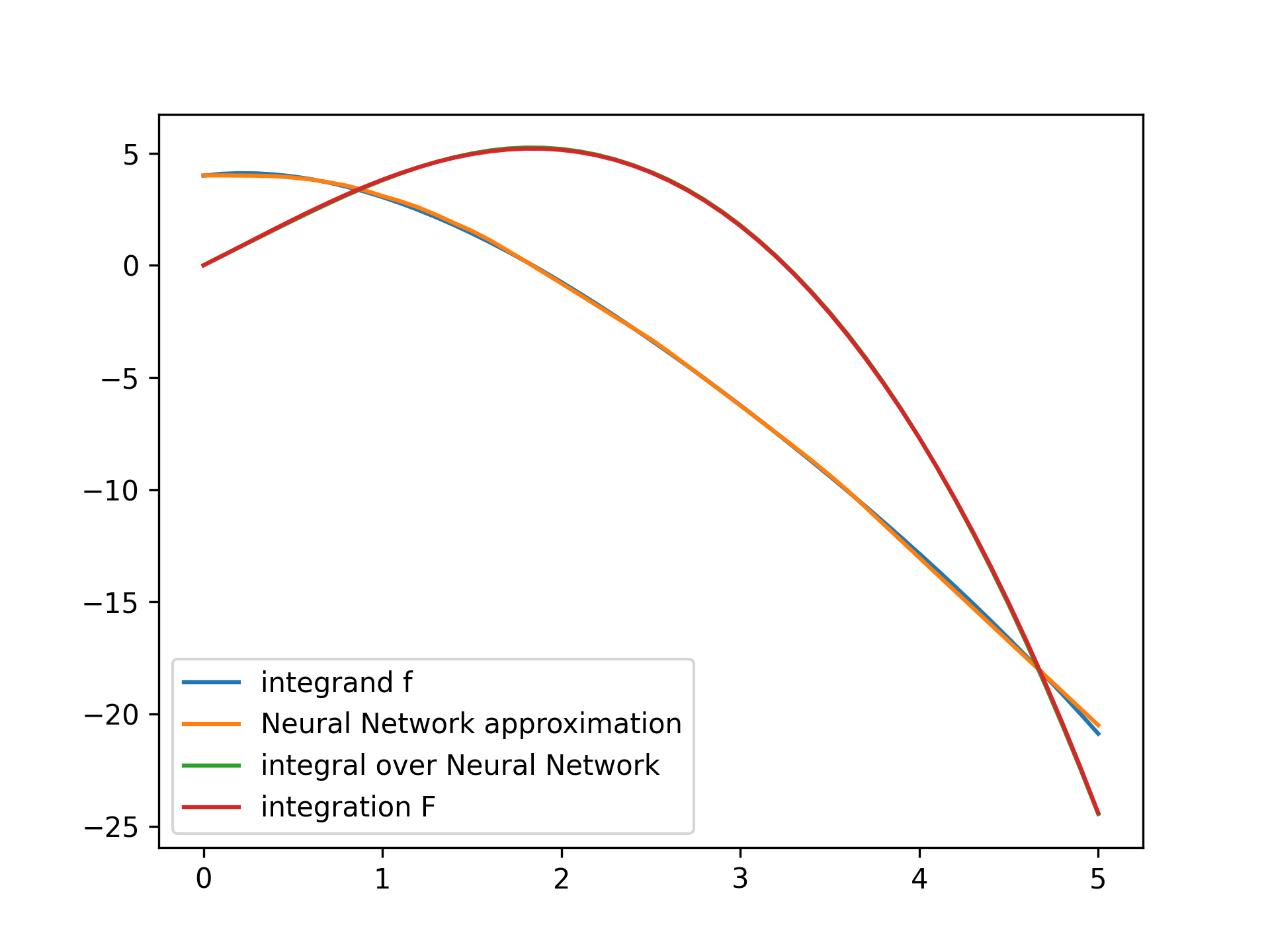}
         \caption{Neural Network approximation}
    \end{subfigure}
    \begin{subfigure}[t]{0.4\textwidth}
         \centering
         \includegraphics[width=\textwidth]{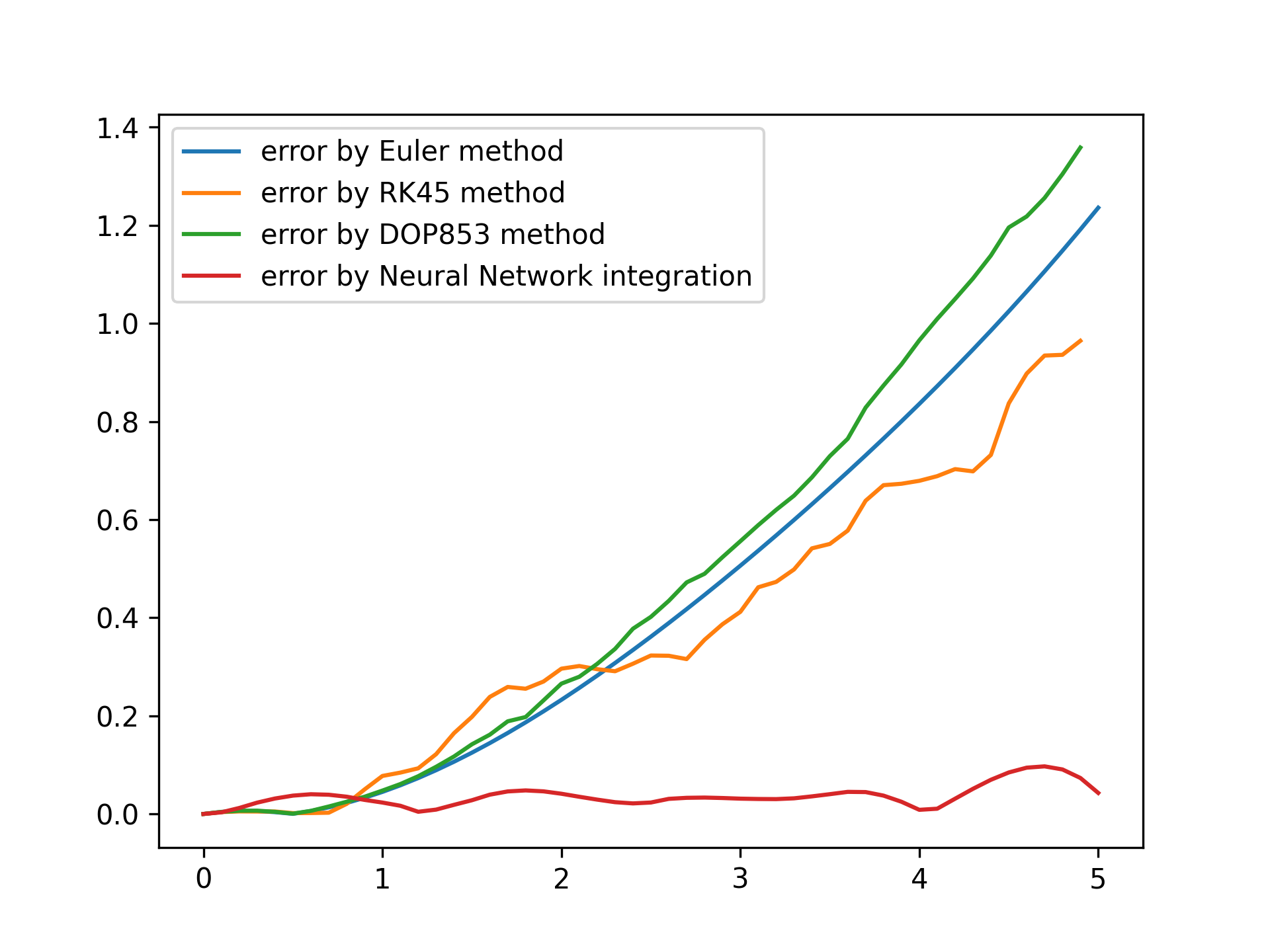}
         \caption{Absolute value of the approximation error given by different numerical integral algorithms.}
    \end{subfigure}
    \caption{Numerical Experiment}
    \label{fig:exp}
\end{wrapfigure}

\section{Integral and Algorithms}
\subsection{Basic Motivation}\label{2}
Let us begin by focusing on a classical numerical question: how can we approximate the integral of an function given only some samples? Assume $f: [a,b] \rightarrow \mathbb{R}$ is a integrable function, denote $F(x) = \int_{a}^{x} f(t)dt$ as the integral of function $f$. Suppose we have $N$ samples of $f$, i.e., $\{\left(x_{n}, f(x_{n}) \right)\}_{n=1}^{N}$. How can we recover $F$ only with these samples without knowing the corresponding values $F(x_{n})$ or the formula of $f$? Especially, can we get an accurate estimation for $F(b) = \int_{a}^{b} f(t)dt$?  

Different from classic numerical integration algorithms, we here introduce one alternative way to solve this problem. With samples, we could approximate $f$ using some integrable estimation $\hat{f}$, then approximate $F$ by integral over $\hat{f}$. Universal Approximation Theorem \citep{cybenko1989approximation,hornik1991approximation,pinkus1999approximation, kidger2020universal} guarantees that, for any $\varepsilon > 0$, there exists a Neural Network $\psi$, such that $|\psi - f| \leq \varepsilon / (b-a)$ on a compact set. Then the integral over $\psi$ satisfies that $|\int_{a}^{x} \psi(t) dt - F(x)| \leq \varepsilon$, which gives a good estimation of function $F$. For high dimensional case, we consider the integral over a closed rectangle.

In the following part, our focus is on obtaining the closed-form solutions for integrating Neural Networks. We first provide explicit integration forms for one layer Neural Networks with any integrable activation function. Then, we derive a piece-wise structure for the integration of multi-layer Neural Networks that use ReLU activation function, and propose an forward integral algorithm with a corrector to enhance the accuracy of the integration.

In Figure \ref{fig:exp}, we illustrate the approximation capabilities of a 2-layer Neural Network, showcasing that our numerical integral approximate perfectly for the integration $F$. Notably, our approach exhibits significantly smaller errors when compared to traditional numerical methods. For details about the experiments, please refer to Appendix \ref{Appd:C}.

\subsection{Integral over one-layer Neural Networks}\label{3.1}
We start with a trivial case: one layer Neural Networks $\psi: \mathbb{R} \rightarrow \mathbb{R}$ 
\begin{equation}
    \begin{aligned}
    y = W^{(1)}x + b^{(1)}, \quad
    \psi(x) = W^{(2)} \sigma(y) + b^{(2)}.
    \end{aligned}
\end{equation}   

Notice that in this one dimensional case, $W^{(2)}$ is a row vector with length $n_{1}$ and $W^{(1)} = (w_{1}, \dots, w_{n_{1}})^\intercal$ is a column vector with length $n_{1}$. Then we state the following Lemma.
\begin{lemma}\label{lemma}
For a one-layer Neural Network $\psi$ defined on a closed interval $[a, b]$, the integral of $\psi$ can be expressed as:
$$
\int_{a}^{x} \psi(t)dt = W^{(2)} z + b^{(2)}(x-a),
$$
where $z = (z_{1}, \dots, z_{n_{1}})^\intercal \in \mathbb{R}^{n_{1}}$ and
$
z_{i} =  \int_{a}^{x} \sigma(w_{i}t + b^{(1)}_{i}) dt.
$
\end{lemma}

\cite{fast} focused on the integral of one layer Neural Network with logistic sigmoid activation function. We demonstrate that our analysis here works for general integrable activation function.

For $n$-dimensional case, the result of integral over rectangles is similar. Since one layer Neural Network is essentially a weighted sum of several integrable function, we only have to repeat Lemma \ref{lemma} for n time. Let $\psi$ be an one-layer Neural Network: $\mathbb{R}^{n} \rightarrow \mathbb{R}$. Then, 
\begin{equation}
    \int_{[a_{1}, b_{1}]\times \dots \times[a_{n}, b_{n}]} \psi(x)dx_{1} \dots dx_{n} 
    =
    \int_{[a_{2}, b_{2}]\times \dots \times[a_{n}, b_{n}]} dx_{2} \dots dx_{n} \{W^{(2)}[\int_{a_{1}}^{b_{1}} \sigma(W^{(1)}x + b^{(1)})dx_{1}] + b^{(2)} (b_{1} - a_{1})\}.
\end{equation}
As long as we have the explicit form of $\tilde{\sigma} = \int \sigma$, we can evaluate the integral directly. For instance, when ignoring the constant term, 
$$
\int \text{ReLU}(x) = \text{ReLU}^{2}(x)/2 \quad \text{and} \quad \int \tanh = \ln{\cosh}.
$$

\subsection{Numerical Integral over multi-layer Neural Networks}\label{3.2}
In this section, we will be delving into the integration of a multi-layer ReLU Neural Network, which is considerably more complex than the one-layer case. A single-layer Neural Network can be expressed as a simple weighted sum of several activation functions, making integration relatively straightforward. However, for a multi-layer Neural Network, the situation is quite different. Even with ReLU activation function, we do not possess explicit knowledge of the areas where the output of a neuron is zero or positive, particularly for neurons in higher layers. As a result, integration of multi-layer Neural Networks is a challenging task.

We start with an observation of ReLU Neural Networks. \citet{arora2018understanding} showed that every ReLU Neural Network is a piece-wise linear function, vice versa. We cite their Theorem here as reference.

\begin{theorem}[\citep{arora2018understanding}]\label{piece}
    Every $\mathbb{R}^{n} \rightarrow \mathbb{R}$ ReLU Neural Network represents a piece-wise linear function, and every piece-wise linear function $\mathbb{R}^{n} \rightarrow \mathbb{R}$ can be represented by a ReLU Neural Network with at most $\lceil\log2(n + 1)\rceil + 1 $ depth.
\end{theorem}

Then to integral a ReLU Neural Network is equivalent with to integral a piece-wise linear function without explicitly knowing the break points. In each piece, the Neural Network can be represented in the form of $\alpha x + \beta$, where $\alpha$ and $\beta$ are coefficients determined by the weights and biases in the network structure. Knowing these coefficients enables us to compute the integral of the Neural Network separately in each piece. Thus, we design a forward integral algorithm. We demonstrate that this algorithm works for batch of input and can be accelerated by GPUs.

\begin{figure}[ht]
    \centering
    \begin{subfigure}[b]{0.3\textwidth}
         \centering
         \includegraphics[width=\textwidth]{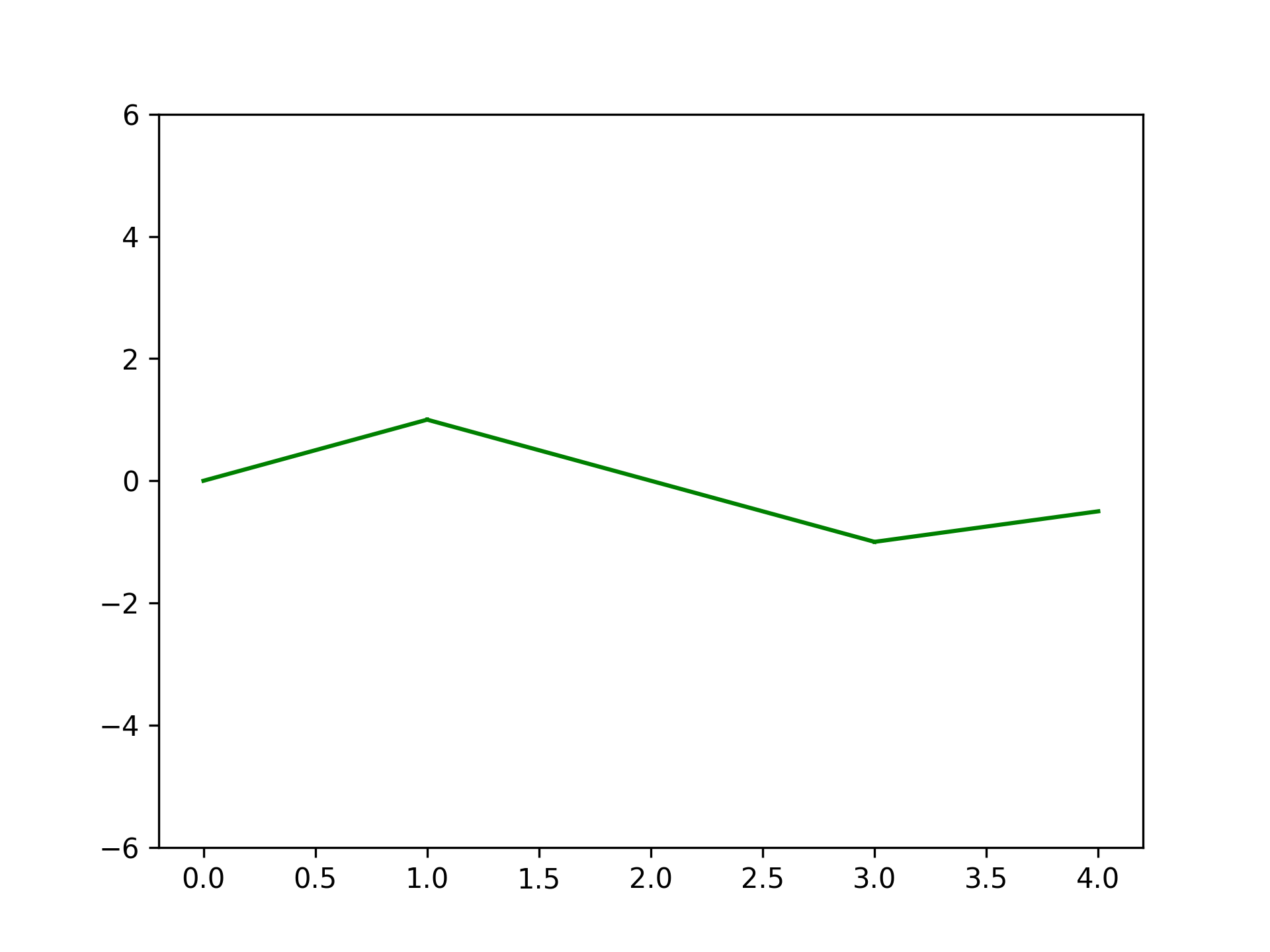}
         \caption{Piece-wise Linear Function}
    \end{subfigure}
    \begin{subfigure}[b]{0.3\textwidth}
         \centering
         \includegraphics[width=\textwidth]{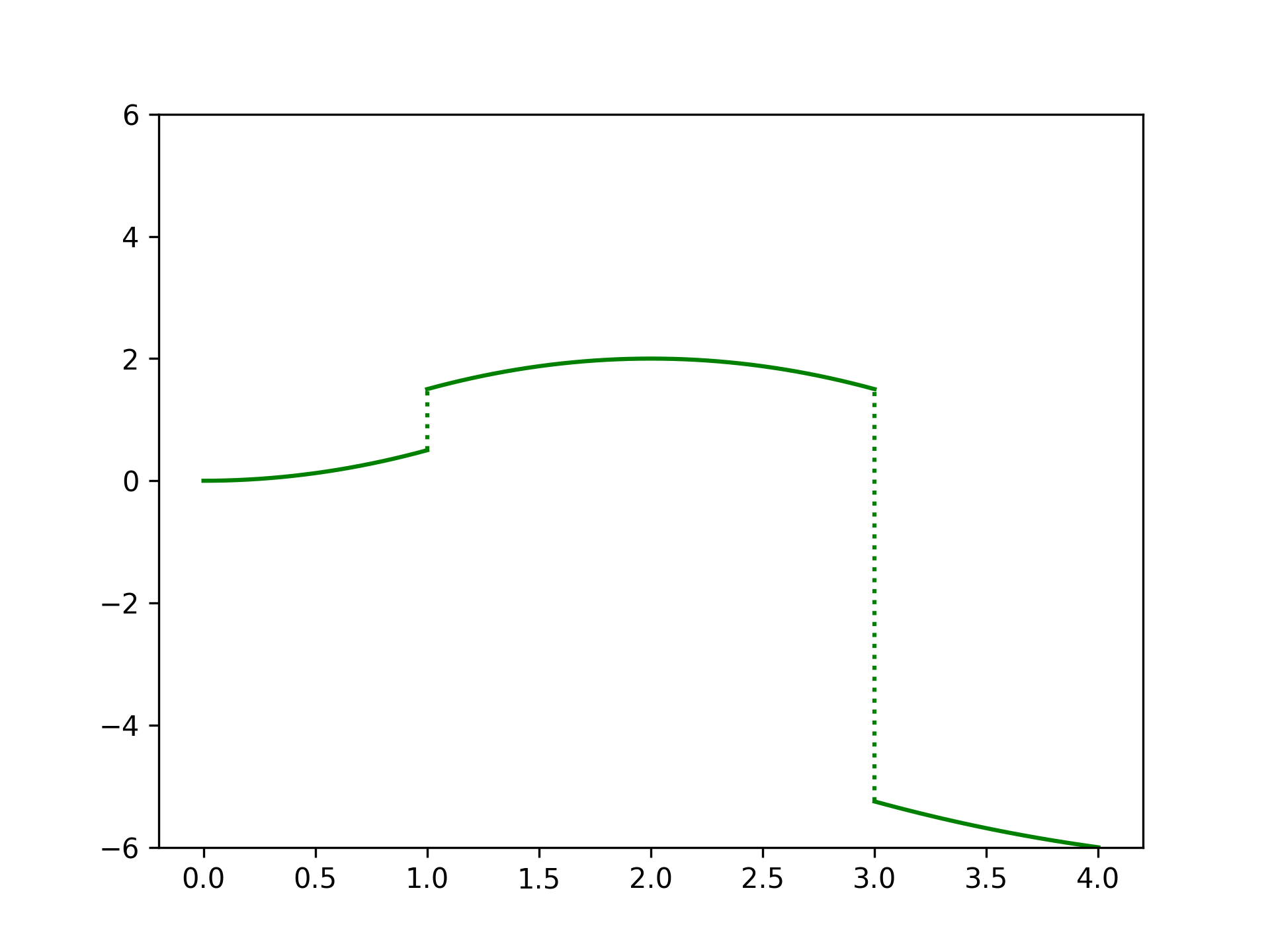}
         \caption{Integral in each piece}
    \end{subfigure}
    \begin{subfigure}[b]{0.3\textwidth}
         \centering
         \includegraphics[width=\textwidth]{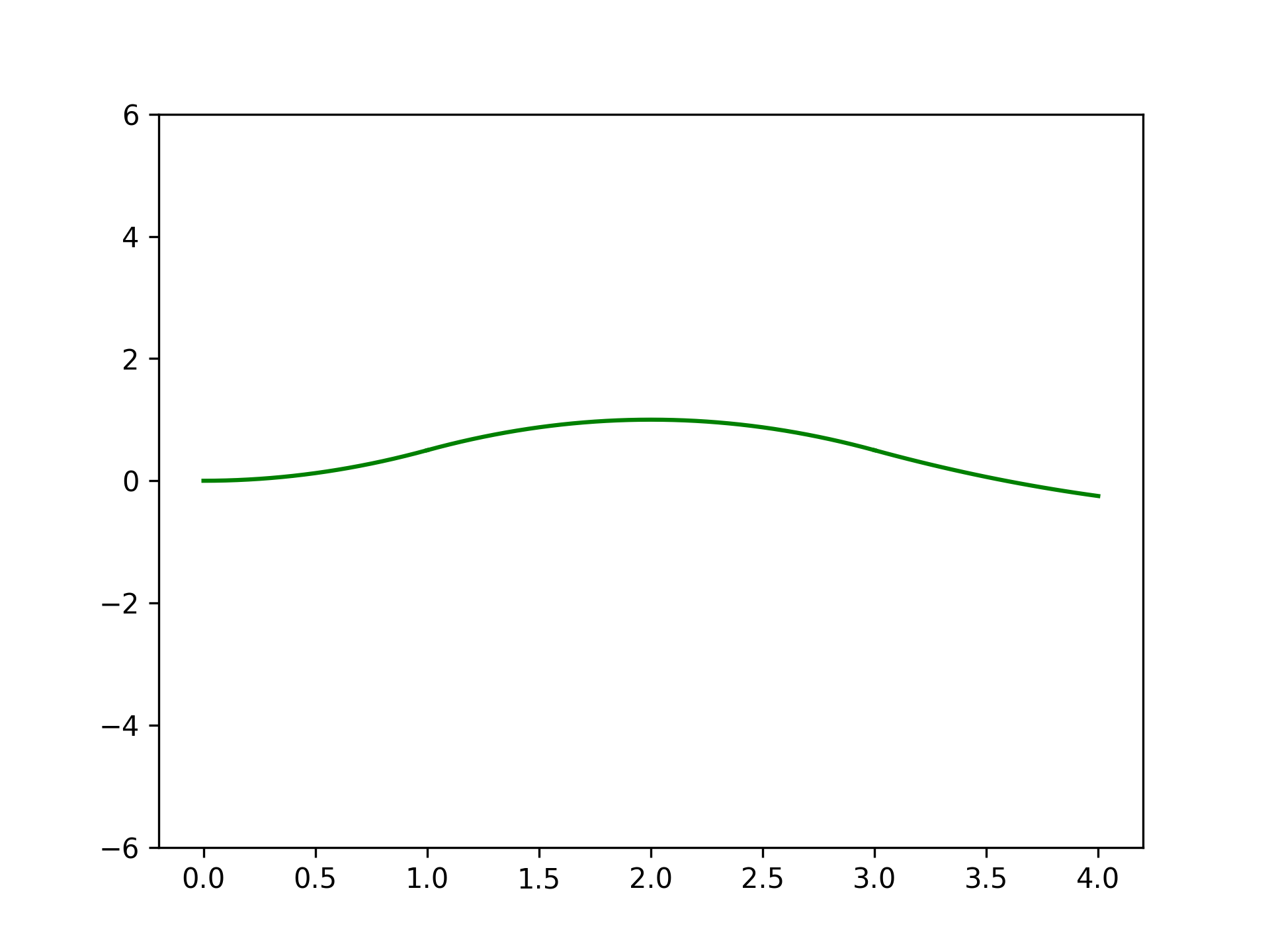}
         \caption{Correct the jump}
    \end{subfigure}
    \caption{Integration Process}
    \label{fig}
\end{figure}

However, by using a piece-wise integral, it is likely that we cannot obtain a continuous function since there are usually jumps at the breakpoints. Nonetheless, it is worth noting that these jumps are actually constant. To address this issue, we propose a numerical algorithm that corrects these jumps and ensures the continuity of the entire integral. Figure \ref{fig} illustrates the entire Neural Network integration process.

\subsubsection{Forward Integral Algorithm}
By Theorem \ref{piece}, when we fix the input $x$, the output of $j$-th neuron in the $i$-th layer could be represented as a linear combination of input $x$ plus a constant term, i.e $\alpha_{ij} x + \beta_{ij}$. Here $\alpha_{ij}$ is a vector and $\beta_{ij}$ is a scalar. Then the output of $(i+1)$-th layer $\alpha_{i+1} x + \beta_{i+1} = \sigma(W^{(i+1)}(\alpha_{i} x + \beta_{i}) + b^{(i+1)})$. Here $\alpha_{i}$ is a matrix and $\beta_{i}$ is a vector. Recall the definition of ReLU, it's equivalent with identity when the input is positive. Otherwise, it always outputs 0. As a result, given input $x$, the coefficient in the next layer depends on the coefficients in the current layer and whether each neuron passes through next layer, i.e
\begin{equation}
    \alpha_{i+1} = \mathbbm{1}_{\{x^{(i+1)}>0\}} \odot (W^{(i+1)} \cdot \alpha_{i}), \quad
    \beta_{i+1} = \mathbbm{1}_{\{x^{(i+1)}>0\}} \odot (W^{(i+1)} \cdot \beta_{i} + b^{(i+1)})  
\end{equation}
Here $\mathbbm{1}$ denotes the indicator function and $\mathbbm{1}_{\{x^{(i)}>0\}}$ represents which neuron passes to next layer, $\odot$ denotes Hadamard product, we use it to ignore neurons which stop. When computing, the Hadamard product between a vector and a matrix, it's implemented for each column of the matrix separately. Thus, following this procedure for each layer, we  summarize our Forward Integral Algorithm \ref{alg:fia} in Appendix \ref{Appd:A}. 

As long as we have $\alpha$ and $\beta$, we can directly get the integral at the piece where input $x$ located in, i.e $\int_{a_{i}}^{x} \psi(x) dx_{i}$ can be represented a Polynomial term plus a constant term. Since the Polynomial term only depends on $\alpha$ and $\beta$, we denote it as Poly[$\alpha$, $\beta$] for convenience. We demonstrate here, this algorithm is just a modified version of the forward algorithm to get output of a Neural Network by adding Hadamard product and indicator operators. 

\subsubsection{Numerical Corrector Algorithm}
By Forward Integral Algorithm, we know the integral at the each piece while ignoring a constant term. Here, we introduce a numerical method to correct the jump at the break point between pieces. 

Given an interval, we first select a partition $\{z_{i}\}_{i=1}^{N}$ of this interval. Denote the corresponding coefficient as $\alpha^{(i)}$ and $\beta^{(i)}$ for each $z_{i}$. Then we connect the integral at $z_{i}$ by adding a constant term $\text{Poly}[\alpha^{(i-1)}, \beta^{(i-1)}](z_{i}) - \text{Poly}[\alpha^{(i)}, \beta^{(i)}](z_{i})$ to make the integral continuous. The corrective steps are outlined in Algorithm \ref{alg:nc} in Appendix \ref{Appd:B}. It's worth noting that $z_{i}$ can be conveniently chosen from samples $x_{i}$.

\subsection{Extension to Convolutional Neural Networks and Residual Neural Networks}\label{3.3}
Convolutional Neural Networks \citep{krizhevsky2017imagenet} and Residual Neural Networks \citep{he2016deep} are more powerful in practice. Our algorithms also work for them, when the activation function is ReLU.

The convolutional layer can be expressed as $\sigma(\kappa * x + b)$, where $*$ denotes the convolution operation, $\kappa$ is the kernel and $b$ is the bias vector. A convolutional Neural Network is one consist of convolutional layers and fully connected layers. 

A Residual Block is usually consists of  convolutional layers, fully-connected layers and a residual connection, i.e $\Psi(x) = x + \psi(x)$, where $\psi(x)$ can be expressed as a convolutional Neural Network. A Residual Neural Network is one taking structures of Residual Blocks, convolutional layers and fully-connected layers. 

We first observe that the convolutional layer and the residual block would not affect the piece-wise structure for ReLU Neural Networks. As a result, our Forward Integral Algorithm works for both ReLU Convolutional Neural Network and ReLU Residual Neural Network. The convolutional layer works as same as the fully connected layer in our Forward Integral Algorithm because the convolution operator by kernel $\kappa$ acts as a weighted linear combination of input. For Residual Neural Networks, we only have to add the coefficients at the beginning of Residual Block and the coefficients of $\psi$ together due to residual connection. Details will be explained in Appendix \ref{Appd:A}.

\section{Discussion}\label{5}
We have developed novel algorithms to compute the integral of Neural Networks. In this section, we outline several promising directions for future research.

One immediate area of interest is the path integrals of Neural Networks, offering potential solutions to challenges in fields such as Molecular Dynamics \citep{marx1996ab,li2022using} and Quantum Mechanics \citep{feynman2010quantum,hj2012introduction}. Our work could be beneficial for these scientific problems, as well as enhancing performance in vision tasks \citep{10.1145/3503250}.

Furthermore, Deep Ritz Method \citep{E2018}, Physics-Informed-Neural-Networks (PINNs) \citep{raissi2019physics,pang2019fpinns,mao2020physics,Karniadakis2021,cai2021physics} and Neural Operators \cite{li2021fourier,Lu2021,kovachki2023neural} have garnered significant attention. While current approaches focus on using derivatives of Neural Networks to satisfy partial differential equations (PDEs), exploring the integration of Neural Networks may open up new avenues for addressing these problems.

Neural Networks have proven effective in addressing challenges associated with density estimation and Bayesian inference, as evidenced by prior research \citep{magdon1998neural, alsing2019fast, lueckmann2019likelihood}. The integrability of Neural Networks holds the promise of substantial improvement in these applications, such as expectation and variance estimation, and entropy estimation, thereby enhancing their statistical perspectives \citep{white1989learning,cheng1994neural}.

Our algorithm depends on how well a Neural Network is trained. So, investigating the error bounds for these Neural Network based numerical algorithm adds an intriguing dimension to our exploration. Understanding the theoretical limits of approximation accuracy \citep{barron1994approximation,ronen2019convergence,devore2021neural} and the factors influencing error is pivotal for assessing the reliability and robustness of such algorithms.

\section*{Acknowledgments}
There is no funding to be disclosed.

\bibliographystyle{plainnat}
\bibliography{references}  

\appendix
\section{Appendix: Forward Integral Algorithm}\label{Appd:A}
We present our Pseudo code for Forward Integral Algorithm \ref{alg:fia} here. Notice that, the last layer is a little different since there is no activation function. Since $x = I \cdot x + \bf{0}$, we initialize $\alpha$ to be identity matrix and $\beta$ to be zero vector $\bf{0}$.

For Residual Neural Network, the algorithm works with little difference. Suppose there is a residual connection after several layers, we can formulate is as $\Psi(x) = \phi(x) + \psi(\phi(x))$, where $\phi$ represents the previous layers, and $\psi$ is a Convolutional Neural Network. Fix a input $x$, denote the coefficients for $\phi$ as $\alpha_{\phi}$, $\beta_{\phi}$ and the coefficients for $\psi(\phi(x))$ as $\alpha_{\psi}$ and $\beta_{\psi}$. Notice that these coefficients can be simply by our Forward Integral Algorithm \ref{alg:fia}. We only have to sum them up. Then the coefficients for $\Psi$ is just $\alpha_{\phi} + \alpha_{\psi}$ and $\beta_{\phi} + \beta_{\psi}$.
\begin{algorithm}[tb]
\caption{Forward Integral Algorithm}\label{alg:fia}
\begin{algorithmic}
\State {\bfseries Input:}  $x \in \mathbb{R}^{n}$, $k$-layer Neural Network with weights $W^{(i)}$, bias $b^{(i)}$
\State {\bfseries Output:}  First order coefficient $\alpha$ and constant term $\beta$
\State $\alpha \gets I$
\State $\beta \gets \bf{0}$
\State $y \gets x$
\For{$i=1$ {\bfseries to} $k+1$}
    \If{$i == k+1$}
        \State $\alpha \gets  W^{(k+1)} \times \alpha$
        \State $\beta \gets W^{(k+1)} \times \beta + b^{(k+1)}$ 
    \Else
        \State $y \gets \sigma(W^{(i)} y + b^{(i)})$
        \State $z \gets \mathbbm{1}_{\{y>0\}}$
        \State $\alpha \gets z \odot (W^{(i)} \times \alpha)$
        \State $\beta \gets z \odot (W^{(i)} \times \beta + b^{(i)})$  
    \EndIf
\EndFor
\State {\bfseries Return:} $\alpha$ and $\beta$
\end{algorithmic}
\end{algorithm}

\section{Appendix: Corrector Algorithm}\label{Appd:B}
Our Corrector Algorithm \ref{alg:nc} lays a crucial role in refining the integral approximation, ensuring continuity at breakpoints within the piece-wise structure. The whole process works as accumulate the constant term through each partition point.
\begin{algorithm}[bt]
\caption{Numerical Corrector}\label{alg:nc}
\begin{algorithmic}
\State {\bfseries Input: Partition $\{z_{i}\}_{i=1}^{N}$ on interval $[a,b]$, k-layer Neural Network $\psi$}
\State {\bfseries Output: $\int_{a}^{b} \psi dx_{j}$}
\State $\alpha_{0} \gets 0$
\State $\beta_{0} \gets 0$
\State $C_{0} \gets 0$
\For{$i=1$ {\bfseries to} $N$}
    \State $\alpha_{i}, \beta_{i} \gets \text{Forward Integral Algorithm}(z_{i})$
    \State $C_{i} \gets  \text{Poly}[\alpha_{i-1}, \beta_{i-1}](z_{i}) - \text{Poly}[\alpha_{i}, \beta_{i}](z_{i}) + C_{i-1}$
\EndFor
\State {\bfseries Return:} $\text{Poly}[\alpha_{N}, \beta_{N}](z_{N}) - \text{Poly}[\alpha_{0}, \beta_{0}](z_{0}) + C_{N}$
\end{algorithmic}
\end{algorithm}

\section{Appendix: Experiment}\label{Appd:C}
We constructed extensive experiments to demonstrate that our algorithms work well on numerical integral, compared with traditional numerical algorithms.

We set the domain to be $[0,5]$, $f(x) = \cos(x) - x^{2} + 4 - 1/(x+1)$. Then the integration $F(x) = \int_{0}^{x} f(t)dt = \sin(x) - 1/3 * x^{3} + 4x - \log(x+1)$.
We trained a 2-layer Neural Network with width $100$ to approximate $f$ with 51 data points $\{x_{i},f(x_{i})\}_{i=0}^{50}$, where $x_{i}$ evenly spaced over $[0, 5]$. We set epochs to be $200$, learning rate $0.001$ with batch size $20$, using Stochastic Gradient Descent and mean squared error loss. When implement our Corrector algorithm, we set partition points $z_{i}$ to be $x_{i}$.

For comparisons, we also implemented Euler-Forward method \citep{euler1845institutionum}, Explicit Runge-Kutta method of order $5(4)$ (RK45) \citep{dormand1980family} and Explicit Runge-Kutta method of order $8$ (DOP853) \citep{wanner1996solving}.

We present our experiment results in Figure \ref{fig:exp}. In the top figure, we display our Neural Network's approximation capabilities. The 2-layer Neural Network exhibits outstanding approximations of the integrand $f$, and its integral simultaneously provides an excellent approximation of $F$. In the bottom figure, we track the approximation error $\left|F(x) - \hat{F}(x)\right|$, the absolute value of difference between true value and estimation. Here $\hat{F}(x)$ represents the estimated integral obtained through numerical algorithms. Notably, the integral computed by the Neural Network displays the smallest approximation error, highlighting its superior performance in our experiment.






\end{document}